\documentclass{article}
\usepackage{amssymb,amsmath}
\usepackage{graphics}

\def\qed{\hfill {\hbox{${\vcenter{\vbox{               
   \hrule height 0.4pt\hbox{\vrule width 0.4pt height 6pt
   \kern5pt\vrule width 0.4pt}\hrule height 0.4pt}}}$}}}

\def\tr{\triangleright}
\def\tl{\triangleleft}

\newtheorem{theorem}{Theorem}
\newtheorem{definition}{Definition}
\newtheorem{lemma}[theorem]{Lemma}
\newtheorem{proposition}[theorem]{Proposition}
\newtheorem{corollary}[theorem]{Corollary}
\newtheorem{example}{Example}
\newtheorem{remark}{Remark}

\newenvironment{proof}[1][Proof]{\smallskip\noindent{\bf #1.}\quad}%
{\qed\par\medskip}

\author{{\begin{tabular}{c} Natasha Harrell \\
\small{\texttt{natasha@math.ucr.edu}}\end{tabular}}
\and
{\begin{tabular}{c} Sam Nelson \\
\small{\texttt{knots@esotericka.org}}\end{tabular}}
\and
\small{\begin{tabular}{c}
Department of Mathematics, University of California, Riverside\\
900 University Avenue, Riverside, CA, 92521\end{tabular} }}

\date{}

\title{\Large \textbf{Quandles and Linking Number}}

\begin{document}
\maketitle

\begin{abstract}
We study the quandle counting invariant for a certain family of 
finite quandles with trivial orbit subquandles. We show how these 
invariants determine the linking number of classical two-component 
links up to sign.
\end{abstract}



\section{\large \textbf{Introduction}}

In the study of finite quandles and the knot invariants they define, 
decomposable quandles are often neglected in favor of indecomposable or 
connected quandles since for knots, two decomposable quandles with the same 
orbit decomposition define the same counting invariant: if $Q=Q_1\cup Q_2$ 
then 
\[\#(\mathrm{Hom}(Q(K),Q)) =\#(\mathrm{Hom}(Q(K),Q_1))
+\#(\mathrm{Hom}(Q(K),Q_2)),\] 
the sum of the counting invariants defined by the orbit subquandles, where
$\#(X)$ denotes the cardinality of the set $X$.

This follows from the fact that knot quandles have a single orbit, so the 
image of any homomorphism $\phi:Q(K)\to T$ from a knot quandle $K$ into a 
finite target quandle $T$ must lie inside a single orbit. The quandle of a 
multi-component link, however, is not connected, and the counting invariants 
defined by nonisomorphic decomposable quandles with the same orbit 
decomposition are generally different for links. 

The motivation for this paper was to study the counting invariants defined 
by decomposable quandles for multi-component virtual links. We started with 
the case of finite target quandles whose orbits are trivial subquandles, so 
that all of the interesting structure is in the off-diagonal portion of the 
quandle's matrix -- trivial subquandles which have been put together in a 
non-trivial way. We call such quandles \textit{trivial orbit quandles} or 
\textit{TOQs}.

Our main result says that the linking number of a two-component classical 
link is determined up to sign by the counting invariants for a family of 
trivial orbit quandles of a particular form.

The paper is organized as follows: in section \ref{q} we review the 
definitions of quandles, the quandle counting invariant and virtual 
linking numbers. In section \ref{l} we give our main result and some 
examples.

\section{\large \textbf{Quandles and the counting invariant}} \label{q}

In this section we review some definitions.

\begin{definition}\textup{
A \textit{virtual link diagram} is a (possibly non-planar) 4-valent
graph with vertices decorated with crossing information, i.e., two
edges meet at a vertex to form an over-strand and the other two form
an under-strand. A \textit{virtual link} is an equivalence class of 
virtual link diagrams under the equivalence relation generated by the
Reidemeister moves (I, II and III). An \textit{oriented} virtual link
diagram has a preferred direction of travel around each component indicated
by arrows. To draw non-planar virtual links on planar paper, we introduce 
\textit{virtual crossings} 
$\left(\raisebox{-0.15in}{\scalebox{0.5}{\includegraphics{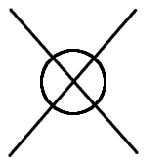}}}
\right)$ to distinguish the crossings arising from non-planarity from the
ordinary \textit{classical crossings}. Two virtual link diagrams are 
equivalent or \textit{virtually isotopic} iff they are related by a finite
sequence of \textit{virtual Reidemeister moves}:
}

\[
\raisebox{-0.2in}{\scalebox{0.6}{\includegraphics{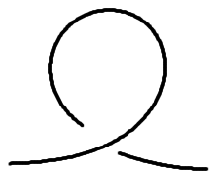}}}
\stackrel{I}{\leftrightarrow}
\raisebox{-0.2in}{\scalebox{0.6}{\includegraphics{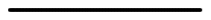}}}
\stackrel{vI}{\leftrightarrow}
\raisebox{-0.2in}{\scalebox{0.5}{\includegraphics{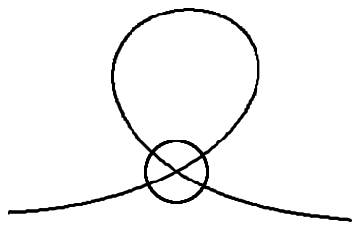}}} \qquad \quad
\raisebox{-0.2in}{\scalebox{0.6}{\includegraphics{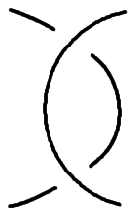}}} 
\stackrel{II}{\leftrightarrow}
\raisebox{-0.2in}{\scalebox{0.6}{\includegraphics{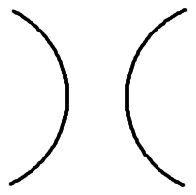}}} 
\stackrel{vII}{\leftrightarrow}
\raisebox{-0.2in}{\scalebox{0.6}{\includegraphics{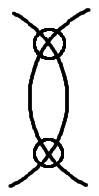}}} \qquad \quad
\]
\[
\raisebox{-0.2in}{\scalebox{0.5}{\includegraphics{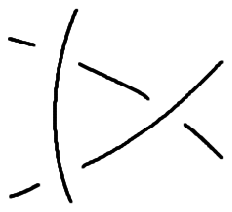}}} 
\stackrel{III}{\leftrightarrow}
\raisebox{-0.2in}{\scalebox{0.5}{\includegraphics{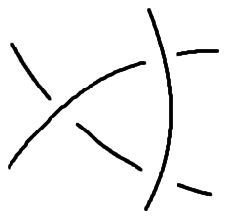}}}\qquad \quad
\raisebox{-0.2in}{\scalebox{0.6}{\includegraphics{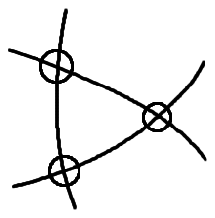}}} 
\stackrel{vIII}{\leftrightarrow}
\raisebox{-0.2in}{\scalebox{0.6}{\includegraphics{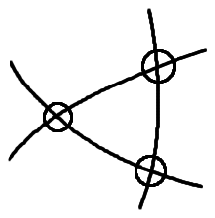}}}\qquad \quad
\raisebox{-0.2in}{\scalebox{0.5}{\includegraphics{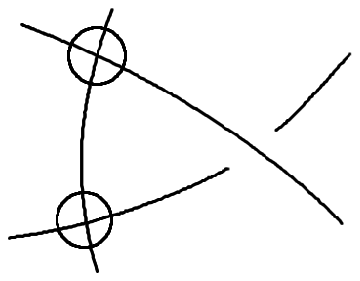}}} 
\stackrel{v}{\leftrightarrow}
\raisebox{-0.2in}{\scalebox{0.5}{\includegraphics{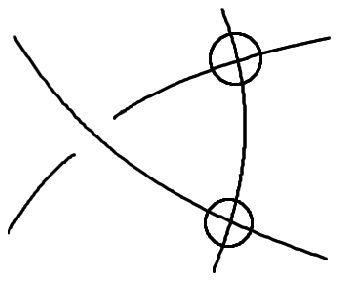}}}
\]

\textup{Classical knot theory may be regarded as a special case of
virtual knot theory, since virtually isotopic classical knots are ambient 
isotopic (see \cite{K}, \cite{GPV}, and \cite{Kb}). In particular, invariants 
of and theorems about virtual knots are invariants of and theorems about 
classical knots.}
\end{definition}

Next, we need

\begin{definition}\textup{
A \textit{quandle} is a set $Q$ endowed with a binary operation 
$\tr:Q\times Q \to Q$ satisfying
\newcounter{qax}
\begin{list}{(\roman{qax})}{\usecounter{qax}}
\item{For all $a\in Q$, $a\tr a = a$,}
\item{For all $a,b\in Q$, there exists a unique $c\in Q$ such that $a=c\tr b$, 
 and}
\item{For all $a,b,c \in Q$, $(a\tr b) \tr c = (a\tr c)\tr (b\tr c)$.}
\end{list}
The uniqueness in axiom (ii) gives us an inverse quandle operation
$\tr^{-1}$ (also denoted by $\bar{\tr}$ or $\tl$), called the \textit{dual}
operation.}

\textup{As expected, a \textit{homomorphism} of quandles is a map 
$f:Q\to Q'$ such that \[f(a\tr b) =f(a)\tr f(b) \quad \forall a,b\in Q.\]
}\end{definition}

We may specify a quandle structure on the finite set $Q=\{x_1,x_2,\dots, x_n\}$
by giving an $n\times n$ matrix $M_Q$ whose entry in row $i$ column $j$ is
$k$ where $x_k=x_i\tr x_j$; that is, $M_Q$ is the operation table of $(Q,\tr)$
where we forget the ``$x$''s and just write the subscripts.

\begin{example}\textup{
The \textit{trivial quandle of order n}, denoted $T_n$, is the quandle with 
matrix
\[
M=\left[
\begin{array}{rrrr}
1 & 1 & \dots & 1 \\
2 & 2 & \dots & 2 \\
\vdots & \vdots & & \vdots \\
n & n & \dots & n  
\end{array}
\right]
\]
}
\end{example}

Quandles and their structure have been studied in many works such as \cite{J}, 
\cite{FR}, \cite{AG} and more. A quandle $Q$ is \textit{decomposable} if we 
may write $Q=Q_1\cup Q_2$ where $Q_1\cap Q_2=\emptyset$ and $Q_1, Q_2$ are 
both subquandles of $Q$. Note that even if a quandle is not decomposable in
this sense, we may still be able to write $Q=\bigcup_{i=1}^n Q_i$ for $Q_i$
disjoint subquandles; see \cite{AG} and \cite{NW}. A quandle is 
\textit{connected} if there is an element $a\in Q$ such that every element of 
$Q$ is equivalent to a word of the form
\[ (\dots ((a \tr^{\pm 1} x_1) \tr^{\pm 1} x_2) \tr^{\pm 1}\dots ) 
\tr^{\pm 1} x_n.\]

In \cite{J} (see also \cite{FR}), a quandle is defined for any oriented tame 
classical link $L$ via a Wirtinger-type presentation, which is called the 
\textit{knot 
quandle} of the link $L$. Specifically, looking in the positive direction of 
the overcrossing strand $b$, if the undercrossing arc on the right is $a$ then 
the one on the left is $a\tr b$. Note that we do not use the
orientation of the undercrossing strand. 
a 

\begin{figure}[!ht] 
\[\scalebox{0.7}{\includegraphics{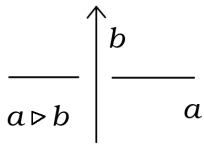}}\]
\caption{Quandle generators at a crossing} \label{crosscond}
\end{figure}

The isomorphism type of the knot quandle is invariant under Reidemeister 
moves; indeed, the quandle axioms are chosen specifically to make this so.  
If $L$ is a knot and $x$ is a generator assigned to an arc in $L$, then 
every generator (and hence every word) in the knot quandle is equivalent 
to a word starting with $x$ since we can just proceed around the knot 
following the specified orientation, picking up quandle operations to 
build a word, until we reach $x$. Hence, the knot quandle of a 
single-component link is connected.

While Joyce defined the knot quandle only for classical knots and links, 
in \cite{K} Kauffman observed that the same Wirtinger-type definition 
extends to virtual knots and links by simply ignoring virtual crossings. 
The resulting quandle, while not a classifying invariant for virtual knots 
and links, is nevertheless an invariant of virtual isotopy.

Given a finite quandle $T$, the number of homomorphisms from a knot quandle
$Q(L)$ into $T$ is an easily computable link invariant: since $Q(L)$ is 
finitely generated, there is a finite set of possible homomorphisms 
into $T$, so we can simply check which of these satisfy the relations in 
$Q(L)$ and hence are homomorphisms. The following is a standard result:

\begin{theorem}
Let $L$ be a link with knot quandle $Q(L)$ and let $T$ be a finite quandle.
Then the number of homomorphisms $\#(\mathrm{Hom}(Q(L),T))$ is an invariant
of quandle isomorphism and hence is an invariant of virtual isotopy of the 
link.
\end{theorem}

Indeed, Joyce \cite{J} showed  that the knot quandle determines both the knot 
group and a longitude, and hence determines the knot type of a classical knot
up to homeomorphism (not necessarily orientation preserving) of 
pairs $(S^3,L)$. Hence, many of the classical link invariants should be 
derivable from the knot quandle. A deeper understanding of how link 
invariants 
arise from the knot quandle should provide us with a better understanding both
of how other link invariants are related and of links themselves. In the next 
section we will show how linking number is determined up to sign by the 
quandle counting invariants associated to a particular kind of non-connected 
quandle.

\section{\large \textbf{Quandle counting and linking number}} \label{l}

We begin with a discussion of linking number for virtual knots.

In \cite{GPV} we find the following definition:

\begin{definition}\textup{
Let $L=L_1\cup L_2$ be an oriented virtual link diagram with two components
$L_1$ and $L_2$. Then $lk_{1/2}(L_1,L_2)$ is the sum of crossing signs at 
crossings where $L_1$ crosses over $L_2$, and $lk_{2/1}(L_1,L_2)$ is the sum 
of crossing signs at crossings where $L_2$ crosses over $L_1$. These are the 
\textit{virtual linking numbers} of $L_1$ with $L_2$. The usual 
(classical) linking number $lk(L_1,L_2) $satisfies
\[lk(L_1,L_2) = \frac{lk_{1/2}(L_1,L_2) +lk_{2/1}(L_1,L_2)}{2}.\]
}\end{definition}


These virtual linking numbers are invariants of virtual isotopy,
since the only virtual move which might potentially change them, move II,
is seen to preserve them by inspection.

For classical links, the two linking numbers are equal since if $L_1$
crosses over $L_2$ at a positive crossing, it must either cross back over 
$L_2$
at a negative crossing (contributing zero to both linking numbers) or
cross under $L_2$ at a positive crossing (contributing $+1$ to both 
linking numbers) to come back to join itself without any virtual crossings. 
For virtual links, the two linking numbers are completely independent. 
Indeed, we have

\begin{proposition}
Let $L=L_1\cup L_2$ be a virtual link. If 
\[lk_{1/2}(L_1,L_2)\ne lk_{2/1}(L_1,L_2)\] then $L$ is non-classical.
\end{proposition}

For virtual links, the virtual linking numbers have some advantages 
over the usual classical linking number. For example, if $L=L_1\cup L_2$ is 
split, then both $lk_{1/2}(L_1,L_2)=0$ and $lk_{2/1}(L_1,L_2)=0$; this is a
stronger condition than $lk(L_1,L_2)=0$, since contributions from two
unequal virtual linking numbers may cancel. 

\begin{example}\textup{
The virtual linking numbers detect the fact that virtual link
\[\includegraphics{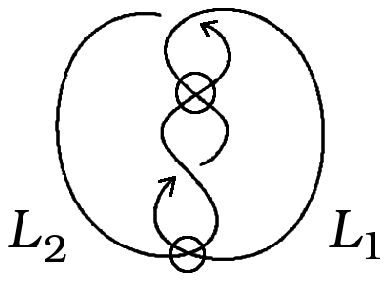}\]
is non-split
since $lk_{1/2}(L_1,L_2)=+1$ while $lk_{2/1}(L_1,L_2)=-1.$ The fact that
$lk_{1/2}(L_1,L_2) \ne lk_{2/1}(L_1,L_2)$ says that $L$ is non-classical.
However, the classical linking number $lk(L_1,L_2)=0.$}
\end{example}

\begin{definition}\textup{
Let $X_n$ be the set $\{1,2,\dots,n+1\}$ with operation $i \tr j=M_{ij}$ 
defined by the matrix}
\[M=\left[
\begin{array}{ccccc}
1 & 1 & \dots & 1 & 2 \\
2 & 2 & \dots & 2 & 3 \\
\vdots & \vdots & & \vdots & \vdots \\
n-1 & n-1 & \dots & n-1 & n \\
n & n & \dots & n & 1 \\
n+1 & n+1 & \dots & n+1 & n+1 
\end{array}
\right].\]
\end{definition}

\begin{proposition}
$X_n$ is a quandle.
\end{proposition}

\begin{remark}
\textup{Note that $X_n$ has two orbit subquandles $S_n\cong T_n$ and 
$S_{n+1}\cong T_1$, both of which are trivial quandles, and we have}
\[ a\tr b = \left\{
\begin{array}{ll}
a+1\ \mathrm{mod} \ n \quad & a\in S_n, \ b= n+1 \\
a & \mathrm{otherwise} 
\end{array}
\right.\] \textup{and}
\[ a\tr^{-1} b = \left\{
\begin{array}{ll}
a-1 \ \mathrm{mod} \ n \quad & a\in S_n, \ b= n+1 \\
a & \mathrm{otherwise}. 
\end{array}
\right.\]
\end{remark}

\begin{proof}

As noted in \cite{HN}, the first two quandle axioms are easy to check since
the diagonal entries are $1,2,3,\dots, n+1$ in order and each column is a 
permutation of $\{1,2,\dots, n+1\}$. Thus, we need only check that the third
quandle axiom is satisfied. We do this in cases. 

\noindent Case 1: $a,b,c$ all in the same subquandle.  In this case,
\[(a\tr b) \tr c = (a\tr c)\tr(b\tr c)\]
since the subquandles are quandles themselves.

\noindent Case 2: $a,b$ in one orbit, $c$ in the other. There are two subcases.
\begin{list}{}{}
\item[(i)]{If $a,b\in S_n$ and $c\in S_{n+1}$ then we have
\[(a\tr b) \tr c = a\tr c = a+1 \ \mathrm{mod}\ n\] and
\[(a\tr c) \tr(b\tr c) = (a+1 \ 
\mathrm{mod} \ n)\tr(b+1 \ \mathrm{mod}\ n ) = a+1 \ \mathrm{mod}\ n.\]}
\item[(ii)]{ If $a,b\in S_{n+1}$ and $c\in S_n$, then
\[(a\tr b) \tr c = a\tr c = a\] and
\[(a\tr c) \tr(b\tr c) = a\tr b = a.\]
}
\end{list}

\noindent Case 3: $b,c$ in one orbit and $a$ in the other. Again we check 
both subcases.
\begin{list}{}{}
\item[(i)]{If $a\in S_n$ and $b,c\in S_{n+1}$ then
\[(a\tr b)\tr c = (a+1 \ \mathrm{mod} \ n) \tr c = a+2 \ \mathrm{mod} \ n
\] and
\[ (a\tr c) \tr(b\tr c)= (a+1 \ \mathrm{mod} \ n) \tr c = a+2 \ \mathrm{mod} \ n.\]
}
\item[(ii)]{If $b,c\in S_n$ and $a\in S_{n+1}$ then
\[(a\tr b) \tr c = a\tr c = a\] and
\[(a\tr c) \tr(b\tr c) = a\tr b = a.\]}
\end{list}

\noindent Case 4: $a,c$ in one orbit and $b$ in the other.
\begin{list}{}{}
\item[(i)]{If $a,c\in S_n$ and $b=S_{n+1}$ then
\[(a\tr b)\tr c = (a+1 \ \mathrm{mod} \ n) \tr c = a+1 \ \mathrm{mod} \ n
\] and
\[ (a\tr c) \tr(b\tr c)= a \tr b = a+1 \ \mathrm{mod} \ n.\]
}
\item[(ii)]{If $a,c\in S_{n+1}$ and $b\in S_n$ then
\[(a\tr b) \tr c = a\tr c = a\] and
\[(a\tr c) \tr(b\tr c) = a\tr (b+1 \ \mathrm{mod} \ n) = a.\]}
\end{list}
\end{proof}

\begin{definition}\textup{
A \textit{coloring} of $L$ by $T$ is a homomorphism of quandles $C:Q(L)\to T$ 
regarded as an assignment of a ``color'' $x\in T$ to each arc in $L$ such that 
the \textit{quandle condition} pictured in figure \ref{crosscond} is satisfied 
at every crossing.
}\end{definition}

\begin{example}\textup{
Every link $L=L_1\cup\dots \cup L_m$ has $n^m$ colorings by $T_n$, the 
trivial quandle of order n. Since $a\tr b=a$ for every $b\in T_n$, all arcs
in the same component must have the same color in any coloring, and we can
choose the colors for each component arbitrarily since every color acts 
trivially on every other color. 
}\end{example}

This example implies that the number of colorings of a link $L$ by
$X_n$ is at least 
\[\#(\mathrm{Hom}(Q(L),S_n))+\#(\mathrm{Hom}(Q(L),S_{n+1})) =n^2+1.\] 
We will now investigate the actual number of colorings of a virtual link 
$L$ by $X_n$.

For the remainder of this paper, let $L$ be a two-component virtual link 
diagram $L=L_1\cup L_2$. To find quandle colorings of $L$, we take a choice of 
color $x\in X_n$ for one arc and attempt to ``push it through'' each crossing,
that is, assign to the opposite undercrossing arc the color determined by the 
quandle coloring condition and the color on the overcrossing arc. Since our 
coloring quandle $X_n$ is a trivial orbit quandle, pushing any color through 
at a crossing where both strands are from the same component yields the same 
color on both underarcs since the colors on a component are from the same 
trivial orbit subquandle. Thus we have:

\begin{lemma}
In a quandle coloring of a link by $X_n$ (or any trivial orbit quandle), the 
colors on the undercrossing strands at a crossing can differ only if the 
overcrossing and undercrossing strands belong to distinct components.
\end{lemma}

The quandle operations in $X_n$ imply the following observation.

\begin{lemma} \label{l0}
In a quandle coloring of $L$ by $X_n$, an inbound undercrossing arc
colored $x\in S_n$ pushes through a crossing with an overarc colored $n+1$ 
to become $x+1\ \mathrm{mod} \ n$ at a positive crossing and 
$x-1\ \mathrm{mod} \ n$ at a negative crossing.
\end{lemma}

We can get an easy upper bound for $\#(\mathrm{Hom}(Q(L),X_n))$:

\begin{proposition}\label{p1}
$L$ has at most $(n+1)^2$ colorings by $X_n$.
\end{proposition}

\begin{proof}
Pushing through any choice of color for a starting arc extends the same color 
along that component of the link until overcrossings from the other component 
of $L$ are reached. For any pair of starting arcs with one arc on each 
component, if both starting colors are from the same orbit, then both colors 
push through unchanged and we have a coloring of the link which is constant 
on each component. If the colors are from different orbits in $X_n$, then one 
component has every arc colored $n+1$, and a given starting coloring 
$a\in S_n$ pushes through an overcrossing colored with $n+1$ to become 
$a+1 \ \mathrm{mod} \ n$ or $a-1 \ \mathrm{mod} \ n$ depending on the 
crossing sign. Continuing pushing through each crossing, we eventually 
reach our original starting arc. Then we have a valid quandle coloring iff 
the color thus determined on the starting arc agrees with the original color. 
In both cases (starting arcs in the same orbit subquandle or in different 
orbit subquandles), the colors on all of the arcs are determined by the 
colors on the initial two arcs. There are $(n+1)^2$ such choices.
\end{proof}

In particular, we have

\begin{lemma} \label{l1}
$L$ has a coloring by $X_n$ with $n+1$ on $L_i$ 
and colors from $S_n$ on $L_j$ iff 
\[lk_{i/j}(L_i,L_j) \equiv 0 \ \mathrm{mod} \ n.\] Moreover, if
$L$ has one such coloring, then it has $n$ such colorings.
\end{lemma}

\begin{proof}
If we color one arc of $L_i$ with $n+1\in X_n$, then in any quandle coloring 
of $L$ every arc of $L_i$ must be colored with $n+1$, since every element 
of $X_n$ acts trivially on $n+1$. 
In order to have a valid quandle coloring, the original color on the 
starting arc must match the color determined by pushing $x$ through the full 
length of the component to return to the starting arc. Lemma \ref{l0} implies 
that this final color is $x+lk_{i/j}(L_i,L_j) \ \mathrm{mod} \ n$, so an 
initial color of $x\ne n+1$ on an arc of $L_j$ pushes 
through to a valid coloring iff $lk_{i/j}(L_i,L_j) \equiv 0 \ \mathrm{mod} 
\ n$. Finally, note that the property of pushing through does not depend on
the actual value of $x\ne n+1$, only its orbit, so if one such color yields
a valid quandle coloring, then all $n$ colors in $S_n$ yield valid quandle 
colorings.
\end{proof}

\begin{corollary}
If $L=L_1\cup L_2$ then \[\#(\mathrm{Hom}(Q(L),X_n)) 
\in\{(n+1)^2, (n+1)^2-n, n^2+1\}.\] 
\end{corollary}

\begin{proof}
We have $\#(\mathrm{Hom}(Q(L),X_n))=(n+1)^2$ if $L$ has colorings with $n+1$
on $L_i$ and colors from $S_n$ on $L_j$ for both $i=1,2$; 
$\#(\mathrm{Hom}(Q(L),X_n))=(n+1)^2-n = n^2+n+1$ if $L$ has colorings with 
$n+1$ on $L_i$ and colors in $S_n$ on $L_j$ for $i=1$ or $i=2$ but not both;
and $\#(\mathrm{Hom}(Q(L),X_n))=n^2+1$ if $L$ has no multi-orbit 
colorings. Note that $\#(\mathrm{Hom}(Q(L),X_n))=(n+1)^2-n $ implies that 
the virtual linking numbers are unequal and $L$ is non-classical.
\end{proof}

\begin{corollary}\label{c1}
$lk_{1/2}(L_1,L_2)=lk_{2/1}(L_1,L_2)=0$ iff
\[\#(\mathrm{Hom}(Q(L),X_n))=(n+1)^2 \quad \mathrm{for \ all} \quad 
n\ge 2.\] 
\end{corollary}

\begin{corollary} \label{c4}
Let $n\in \mathbb{Z}$ with $n\ge2.$ If 
\[|lk_{1/2}(L_1,L_2)|+ |lk_{2/1}(L_1,L_2)|\ne 0\] 
then
\[\#(\mathrm{Hom}(Q(L),X_n))=(n+1)^2 \iff n|\ |lk_{1/2}(L_1,L_2)| 
\ \mathrm{and} \ n|\ |lk_{2/1}(L_1,L_2)|.\]
\end{corollary}

\begin{proof}
$\#(\mathrm{Hom}(Q(L),X_n))=(n+1)^2$ if and only if each starting coloring
choice for a given pair of arcs on the two components pushes through
to a valid coloring. By lemma \ref{l1}, this can only happen if both
virtual linking numbers are zero mod $n$, i. e., if 
$n|\ |lk_{1/2}(L_1,L_2)|$ 
and $\ n|\ |lk_{2/1}(L_1,L_2)|.$
\end{proof}

\begin{corollary}\label{c3}
If \[|lk_{1/2}(L_1,L_2)|+|lk_{2/1}(L_1,L_2)|\ne 0\] then 
$\#(\mathrm{Hom}(Q(L),X_n))<(n+1)^2$ for any $n$ greater than the number
of two-component crossings in the link diagram $L$.
\end{corollary}

\begin{corollary} \label{c2}
If either $|lk_{1/2}(L_1,L_2)|=1$ or $|lk_{2/1}(L_1,L_2)|=1$
then $\#(\mathrm{Hom}(Q(L),X_n))<(n+1)^2$ for all $n\ge 2.$ Conversely,
if $\#(\mathrm{Hom}(Q(L),X_n))<(n+1)^2$ for all $n\ge 2$ then either
$|lk_{1/2}(L_1,L_2)|=1$, $|lk_{2/1}(L_1,L_2)|=1$ or both 
$|lk_{1/2}(L_1,L_2)|=|lk_{2/1}(L_1,L_2)|=1$.
\end{corollary}

\begin{proof}
If $|lk_{i/j}(L_1,L_2)|=1$ then 
$|lk_{1/2}(L_1,L_2)|+|lk_{1/2}(L_1,L_2)|\ne 0$. Then  $n\ge 2$ in
corollary \ref{c4} implies that $\#(\mathrm{Hom}(Q(L),X_n))<(n+1)^2$,
since $n\not| 1$.

Conversely, suppose $\#(\mathrm{Hom}(Q(L),X_n))<(n+1)^2$ for all $n\ge 2$.
Then by corollary \ref{c1} at least one of the two virtual linking
numbers is greater than 0 in absolute value. If both are greater than
1 in absolute value, then let $N=\gcd(|lk_{1/2}(L_1,L_2)|, 
|lk_{2/1}(L_1,L_2)|)$. Then by lemma \ref{l1}, if $N>1$ then 
$\#(\mathrm{Hom}(Q(L),X_N))=(N+1)^2$, contrary to hypothesis. Hence
either $|lk_{1/2}(L_1,L_2)|=1$, $|lk_{2/1}(L_1,L_2)|=1$, or both
$|lk_{1/2}(L_1,L_2)|=|lk_{2/1}(L_1,L_2)|=1$.
\end{proof}

\begin{theorem}
Let $L=L_1\cup L_2$ be a virtual link diagram and let 
\begin{eqnarray*}
S & = & \{n\ : \ \#(\mathrm{Hom}(Q(L),X_n))=(n+1)^2, n\ge 2\} \ 
\mathrm{and} \\
S' & = &  \{n\ : \ \#(\mathrm{Hom}(Q(L),X_n))=(n+1)^2 -n, n\ge 2\}.
\end{eqnarray*}

Then if $0<\#(S)<\infty$ we have
\[\max(S)=\mathrm{gcd}(|lk_{1/2}(L_1,L_2)|,|lk_{2/1}(L_1,L_2)|)\] and
\[\max(S')=\max(\{ n\ :\ n| \ |lk_{i/j}(L_1,L_2)| \ \mathrm{and} 
\ n\not|  \ |lk_{j/i}(L_1,L_2)|\}).\]
\end{theorem}

\begin{proof}
If $0<\#(S)<\infty$ then neither linking number is 0 or 1 in absolute value 
by corollaries
\ref{c1} and \ref{c2}. Then $n\in S'$ says that $n$ divides one linking 
number but not the other by lemma \ref{l1}; $\max(S)$ is the largest
number which divides both $|lk_{1/2}(L_1,L_2)|$ and $|lk_{2/1}(L_1,L_2)|$
by corollary \ref{c4}.
\end{proof}

We can now state our main result, which says that the classical linking
number is determined up to sign by the quandle coloring invariants with
target quandle $X_n$ for $n\ge 2$.

\begin{theorem} \label{main}
Let $L=L_1\cup L_2$ be a classical link diagram and let 
\[S=\{n\ : \ \#(\mathrm{Hom}(Q(L),X_n))=(n+1)^2,\  n\ge 2\}.\] Then 
\[|lk(L_1,L_2)| = \frac{|lk_{1/2}(L_1,L_2)|+|lk_{2/1}(L_1,L_2)|}{2}
= \left\{ \begin{array}{ll}
0, & \#(S)=\infty, \\
1, & \#(S)=0,\  \mathrm{or} \\
\max(S), & 0<\#(S)<\infty
\end{array} \right.
\]
\end{theorem}

\begin{proof}
If both virtual linking numbers are zero, then corollary \ref{c1} says
$S=\{n\in \mathbb{Z} \ | \ n\ge 2\}$. If both are $\pm 1$, then by corollary 
\ref{c2}
we have $S=\emptyset$. Finally, if $|lk_{1/2}(L_1,L_2)|=m$ (which equals
$|lk_{2/1}(L_1,L_2)|$ since $L$ is classical) then the set $S$ includes all 
and only positive integer divisors of $m$, so $m =\max (S)$. 
\end{proof}

\begin{example}
\[\raisebox{-0.9in}{\includegraphics{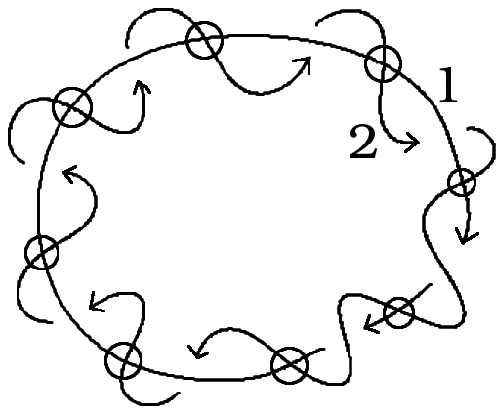}}  \quad \quad
\begin{array}{|r|c|} \hline
n & \mathrm{Hom}(Q(L),X_n)  \\ \hline
2 & \ 9 = (n+1)^2 \\
3 & 13 = (n+1)^2-n \\
4 & 17 = (n+1)^2-2n \\
5 & 26 = (n+1)^2-2n \\
6 & 43 = (n+1)^2-n \\ 
7 & 50 = (n+1)^2-2n  \\ \hline
\end{array}
\]

\textup{This virtual link has two components labeled $1$ and $2$, with virtual
linking numbers $lk_{1/2}(1,2)=6$ and $lk_{2,1}(1,2)=-2.$ Then 
$\mathrm{gcd}(|lk_{1/2}(L_1,L_2)|, |lk_{2/1}(L_1,L_2)|) = 2$ and
the largest cardinality which divides one linking number but not both is
$6$. Our \textit{Maple} computations confirm the values in the table.}

\end{example}

As a practical observation, we note that corollary \ref{c3} implies that it is
not necessary to test an infinite number of quandle counting invariants to 
determine $\#(S)$; it is sufficient to compute $\#(\mathrm{Hom}(Q(L),X_n))$ 
for $2\le n\le N$ where $N$ is the number of multi-component crossings in $L$.
Of course, our interest in theorem \ref{main} lies not in finding a new
less practical method of computing linking numbers but rather as a connection
between quandle-based link invariants and classical link invariants.


\begin{thebibliography}{00}

\bibitem{AG}{N. Andruskiewitsch and M. Gra\~{n}a.
 From racks to pointed Hopf algebras.
 \textit{Adv. Math.} \textbf{178} (2003) 177-243.}

\bibitem{FR}{ R. Fenn and C. Rourke. Racks and links in codimension two.
\textit{ J. Knot Theory Ramifications}  \textbf{1}  (1992) 343-406.}

\bibitem{HN}{B. Ho and S. Nelson. Matrices and finite quandles. 
\textit{Homol. Homotop. Appl.} \textbf{7} (2005) 197-208.}

\bibitem{GPV}{ M. Goussarov, M. Polyak and O. Viro. 
Finite-type invariants of classical and virtual knots. \textit{Topology} 
\textbf{39}  (2000) 1045-1068. }

\bibitem{J}{D. Joyce.
 A classifying invariant of knots, the knot quandle.
 \textit{J. Pure Appl. Algebra}  \textbf{23}  (1982)  37-65.}

\bibitem{K}{L. H. Kauffman. Virtual Knot Theory. 
\textit{ European J. Combin.} \textbf{20}  (1999) 663-690.}

\bibitem{Kb}{ G. Kuperberg. What is a virtual link? \textit{Algebr. 
Geom. Topol.} \textbf{3} (2003) 587-591.}

\bibitem{NW}{S. Nelson and C-Y. Wong. On the orbit decomposition of finite
quandles. \textit{J. Knot Theory Ramifications} \textbf{15} (2006) 761-772.}

\end{thebibliography}
\end{document}